\documentclass[12pt]{amsart}
\usepackage{amsmath,amsfonts,amssymb,amscd,mathrsfs,pstricks,pst-plot}

\input xy
\xyoption{all}

\newtheorem{theorem}{Theorem}[section]

\newtheorem{corollary}[theorem]{Corollary}

\theoremstyle{definition}
\newtheorem{definition}[theorem]{Definition}
\newtheorem{example}[theorem]{Example}

\numberwithin{equation}{section}

\newcommand{\C}{\mathbb{C}}
\renewcommand{\P}{\mathbb{P}}
\newcommand{\R}{\mathbb{R}}
\newcommand{\sO}{\mathscr{O}}
\newcommand{\sC}{\mathscr{C}}
\newcommand{\sM}{\mathscr{M}}
\newcommand{\sS}{\mathscr{S}}

\begin{document}

\title{Survey of Oka theory}

\author{Franc Forstneri\v c and Finnur L\'arusson}

\address{Franc Forstneri\v c, Faculty of Mathematics and Physics, University of Ljubljana, and Institute of Mathematics, Physics and Mechanics, Jadranska 19, 1000 Ljubljana, Slovenia}
\email{franc.forstneric@fmf.uni-lj.si}

\address{Finnur L\'arusson, School of Mathematical Sciences, University of Adelaide, Adelaide SA 5005, Australia}
\email{finnur.larusson@adelaide.edu.au}

\thanks{The first-named author was partially supported by grants P1-0291 and J1-2043-0101, ARRS, Republic of Slovenia, and by the conference RAFROT, Rinc\'on, Puerto Rico, March 2010.}


\date{10 September 2010.  Last minor change 19 December 2010.}


\begin{abstract}
Oka theory has its roots in the classical Oka principle in complex analysis.  It has emerged as a subfield of complex geometry in its own right since the appearance of a seminal paper of M.\ Gromov in 1989.  Following a brief review of Stein manifolds, we discuss the recently introduced category of Oka manifolds and Oka maps.  We consider geometric sufficient conditions for being Oka, the most important of which is ellipticity, introduced by Gromov.  We explain how Oka manifolds and maps naturally fit into an abstract homotopy-theoretic framework.  We describe recent applications and some key open problems.  This article is a much expanded version of the lecture given by the first-named author at the conference RAFROT 2010 in Rinc\'on, Puerto Rico, on 22 March 2010, and of a recent survey article by the second-named author \cite{Larusson6}.
\end{abstract}

\maketitle
\tableofcontents


\section{Introduction}
\label{introduction}

\noindent
The Oka principle first appeared in the pioneering work of K.\ Oka, who proved in 1939 that a second Cousin problem on a domain of holomorphy can be solved by holomorphic functions if it can be solved by continuous functions \cite{Oka1939}.  It follows that in a holomorphic fibre bundle with fibre $\C^*$ over a domain of holomorphy, every continuous section can be continuously deformed to a holomorphic section.  In the following decades further results in a similar spirit were proved, notably by H.\ Grauert, who extended Oka's theorem to more general fibre bundles over Stein spaces \cite{Grauert3}.  The idea emerged of a general Oka principle, saying that on Stein spaces, there are only topological obstructions to solving holomorphic problems that can be cohomologically formulated.

The modern development of the Oka principle started with M.\ Gromov's seminal paper of 1989  \cite{Gromov:OP}, in which the emphasis moved from the cohomological to the homotopy-theoretic aspect.  As a result, Oka's theorem was extended once again to more general fibre bundles and even to certain maps that are not locally trivial.  A major application, the solution of Forster's conjecture, appeared a few years later.  Since 2000, there has been an ongoing effort to further develop the theory and applications of Gromov's Oka principle.

Our goal is to survey the present state of Oka theory.  Here is a brief summary of where we stand as of August 2010, to be fleshed out in detail in the rest of the paper.

\noindent
$\bullet$ \ There is a good definition of a complex manifold being Oka.  This notion is, in a sense, dual to being Stein and opposite to being Kobayashi hyperbolic.  It naturally extends to holomorphic maps and has many nontrivially equivalent formulations.  Sections of Oka maps satisfy a strong Oka principle.

\noindent
$\bullet$ \ In essence, Gromov's Oka principle is about sufficient geometric conditions for the Oka property to hold.  Gromov's results in this direction have been strengthened considerably.

\noindent
$\bullet$ \ There is a growing list of applications of Gromov's Oka principle for which the older theory does not suffice.

\noindent
$\bullet$ \ There is an underlying model structure that shows that the Oka property is, in a precise sense, homotopy-theoretic.

All this appears to us evidence that Oka theory has reached a certain maturity and will play a future role as a subfield of complex geometry in its own right.  Many problems remain open; the ones that seem most important to us are listed at the end of the paper.


\section{Stein manifolds}
\label{Stein}

\noindent
The basic theory of Stein manifolds was developed in the period 1950 to 1965.  We survey the highlights and mention a few more recent developments that are relevant here.

The central concept of classical several complex variables is that of a domain of holomorphy: a domain in complex affine space $\C^n$ with a holomorphic function that does not extend holomorphically to any larger domain, even as a multivalued function.  Every domain in $\C$ is a domain of holomorphy.  One of the discoveries that got several complex variables started at the turn of the 20th century is that this is far from true in higher dimensions.  

The notion of a Stein manifold, introduced by K.\ Stein in 1951, generalizes domains of holomorphy to the setting of complex manifolds (here always assumed second countable, but not necessarily connected).  Roughly speaking, Stein manifolds are the complex manifolds whose function theory is similar to that of domains in $\C$.  There are at least four fundamentally different ways to precisely define the concept of a Stein manifold.  The equivalence of any two of these definitions is a nontrivial theorem.

First, Stein's original definition, simplified by later developments, states that a complex manifold $S$ is {\it Stein} (or, as Stein put it, {\it holomorphically complete}) if it satisfies the following two conditions.
\begin{enumerate}
\item  Holomorphic functions on $S$ separate points, that is, if $x,y\in S$, $x\neq y$, then there is $f\in\sO(S)$ such that $f(x)\neq f(y)$.  Here, $\sO(S)$ denotes the algebra of holomorphic functions on $S$.
\item  $S$ is {\it holomorphically convex}, that is, if $K\subset S$ is compact, then its $\sO(S)$-hull $\widehat K_{\sO(S)}$, consisting of all $x\in S$ with $\left|f(x)\right|\leq\max_K\left|f\right|$ for all $f\in\sO(S)$, is also compact.  Equivalently, if $E\subset S$ is not relatively compact, then there is $f\in\sO(S)$ such that $f|E$ is unbounded.
\end{enumerate}
A domain in $\C^n$ is Stein if and only if it is a domain of holomorphy.  Every noncompact Riemann surface is Stein.  

Second, a connected complex manifold is Stein if and only if it is biholomorphic to a closed complex submanifold of $\C^m$ for some $m$.  Namely, submanifolds of $\C^m$ are clearly Stein: the functions $f$ in the definition above can be taken to be coordinate projections.  (More generally, it is easy to see that a closed complex submanifold of a Stein manifold is itself Stein.)  Conversely, R.\ Remmert proved in  1956 that every connected Stein manifold $S$ admits a proper holomorphic embedding into $\C^m$ for some $m$.  In 1960--61, E.\ Bishop and R.\ Narasimhan independently showed that if $\dim_\C S=n$, then $m$ can be taken to be $2n+1$.  The optimal embedding result is that if $n\geq 2$, then $m$ can be taken to be $\left[3n/2\right]+1$.  This was a conjecture of O.\ Forster, who showed that for each $n$, no smaller value of $m$ works in general.  Forster's conjecture was proved in the early 1990s by Y.\ Eliashberg and M.\ Gromov \cite{Eliashberg-Gromov2} (following their much earlier paper \cite{Eliashberg-Gromov1}) and J.\ Sch\"urmann \cite{Schurmann}.  The proof relies on Gromov's Oka principle discussed below.  This problem is still open in dimension 1.  Every noncompact Riemann surface properly embeds into $\C^3$, but relatively few are known to embed, even non-properly, into $\C^2$.  For recent results in this direction, see \cite{Forstneric-Wold1, Majcen-2009}.

Third, Stein manifolds are characterized by a cohomology vanishing property.  The famous Theorem B of H.\ Cartan, proved in his seminar in 1951--52 as the first triumph of the sheaf-theoretic approach to complex analysis, states that if a complex manifold $S$ is Stein, then $H^k(S,\mathscr F)=0$ for every coherent analytic sheaf $\mathscr F$ on $S$ and every $k\geq 1$.  The converse is easy.

Finally, Stein manifolds can be defined in terms of plurisubharmonicity.  This notion, which is ordinary convexity (in some holomorphic coordinates) weakened just enough to be preserved by biholomorphisms, was introduced independently by P.\ Lelong and by K.\ Oka in 1942 and plays a diverse and fundamental role in higher-dimensional complex analysis.  H.~Grauert proved in 1958 that a complex manifold $S$ is Stein if and only if there is a strictly plurisubharmonic function $\rho:S\to\R$ (smooth if desired, or just upper semicontinuous) that is an exhaustion in the sense that for every $c\in \R$, the sublevel set $\{x\in S : \rho(x) < c\}$ is relatively compact in $S$.  This result is the solution of the Levi problem for manifolds.

If $S$ is a Stein manifold embedded as a closed complex submanifold of $\C^m$, then the restriction to $S$ of the square $\left\Vert\boldsymbol\cdot\right\Vert^2$ of the Euclidean norm is a smooth strictly plurisub\-harmonic exhaustion (and the translate $\left\Vert\boldsymbol\cdot-a\right\Vert^2$, for a generic point $a\in\C^m$, is in addition a Morse function).  The converse implication uses a deep result, proved by Grauert and, using Hilbert space methods, by L.\ H\"ormander in 1965, that the existence of a strictly plurisubharmonic exhaustion $\rho$ on a complex manifold $S$ implies the solvability of all consistent $\overline\partial$-equations on $S$.  In particular, $\rho^{-1}(-\infty,c_1)$ is Runge in $\rho^{-1}(-\infty,c_2)$ for all real numbers $c_1<c_2$, from which the defining properties (1) and (2) of a Stein manifold easily follow.

The definition of a complex space (not necessarily reduced) being Stein is the same as for manifolds.  Theorem B still holds.  The property of having a continuous strictly plurisubharmonic exhaustion is still equivalent to being Stein (Narasimhan 1962).  And a complex space is Stein and has a bound on the dimensions of its tangent spaces if and only if it is biholomorphic to a complex subspace of $\C^m$ for some $m$ (Narasimhan 1960).

Often the best way to show that a complex space is Stein is to produce a strictly plurisubharmonic exhaustion on it.  For example, this is how Y.-T.\ Siu proved in 1976 that a Stein subvariety of a reduced complex space has a basis of Stein open neighbourhoods \cite{Siu1976}.  Stein neighbourhood constructions sometimes allow us to transfer a problem on a complex manifold to Euclidean space where it becomes tractable.  A recent example is the application of the Stein neighbourhood construction in \cite{Forstneric-Wold2} to the proof that the basic Oka property implies the parametric Oka property for manifolds \cite{FF:OkaManifolds}.

By Morse theory, the existence of a smooth strictly plurisubharmonic Morse exhaustion on a Stein manifold $S$ of complex dimension $n$ has the important topological consequence that $S$ has the homotopy type of a CW complex of real dimension at most $n$.  The reason is that the Morse index of a nondegenerate critical point of a plurisubharmonic function on $S$ is at most $n$.  This simple observation has a highly nontrivial converse, Eliashberg's topological characterization of Stein manifolds of dimension at least $3$ \cite{Eliashberg1}:  If $(S,J)$ is an almost complex manifold of complex dimension $n\geq 3$, which admits a Morse exhaustion function $\rho:S\to\R$ all of whose Morse indices are at most $n$, then $J$ is homotopic to an integrable complex structure $\tilde J$ on $S$ in which $\rho$ is strictly plurisubharmonic, so $(S,\tilde J)$ is Stein.  

Eliashberg's result fails in dimension $2$.  The simplest counterexample is the smooth manifold $S^2\times\R^2$, which satisfies the hypotheses of the theorem but carries no Stein structure by Seiberg-Witten theory.  Namely, $S^2\times\R^2$ contains embedded homologically nontrivial spheres $S^2\times\{c\}$, $c\in\R^2$, with self-intersection number $0$, while the adjunction inequality shows that any homologically nontrivial sphere $C$ in a Stein surface has $C\cdot C\leq -2$.  However, Gompf has shown that Eliashberg's theorem holds up to homeomorphism in dimension $2$ \cite{Gompf1,Gompf2}.  More precisely, a topological 4-manifold is homeomorphic to a Stein surface if and only if it is oriented and is the interior of a (possibly infinite) topological handlebody with only $0$-, $1$-, and $2$-handles.  It follows that there are Stein surfaces homeomorphic to $S^2\times\R^2$.

The monographs by Grauert and Remmert \cite{Grauert-Remmert}, by Gunning and Rossi \cite{Gunning-Rossi}, and by H\"ormander \cite{Hormander-SCV}, now considered classics, are still excellent sources for the theory of Stein manifolds and Stein spaces.


\section{Oka manifolds}
\label{Oka}

\noindent
Roughly speaking, Stein manifolds are characterized by carrying many holomorphic functions, that is, by having many holomorphic maps into $\C$.  Dually, Oka manifolds have many holomorphic maps from $\C$, and more generally from Stein manifolds, into them.  To introduce the precise definition of an Oka manifold, we begin with two well known theorems of 19th century complex analysis.

\noindent{\bf Runge Theorem.}  If $K$ is a compact subset of $\C$ with no holes, that is, the complement $\C\setminus K$ is connected, then every holomorphic function on $K$ can be approximated uniformly on $K$ by entire functions.  (By a holomorphic function on $K$ we mean a holomorphic function on some open neighbourhood of $K$.)

\noindent{\bf Weierstrass Theorem.}  On a discrete subset of a domain $\Omega$ 
in $\C$, we can prescribe the values of a holomorphic function on $\Omega$.

In the mid-20th century, these results were generalized to Stein manifolds.

\noindent{\bf Oka-Weil Approximation Theorem.}  If a compact subset $K$ of a Stein manifold $S$ is holomorphically convex, that is, $K=\widehat K_{\sO(S)}$, then every holomorphic function on $K$ can be approximated uniformly on $K$ by holomorphic functions on $S$.

A compact subset $K$ of $\C$ is holomorphically convex if and only if it has no holes.  In higher dimensions, holomorphic convexity is much more subtle; in particular, it is not a topological property.

\noindent{\bf Cartan Extension Theorem.}  If $T$ is a closed complex subvariety of a Stein manifold $S$, then every holomorphic function on $T$ extends to a holomorphic function on $S$.

We usually consider these theorems as fundamental results about Stein manifolds, but they can also be viewed as expressing properties of the target manifold, the complex number field $\C$.  Replacing $\C$ by an arbitrary complex manifold $X$, we are led to the following two properties that $X$ may or may not have.  To avoid topological obstructions, which are not relevant here, we restrict ourselves to very special $S$, $K$, and $T$.

\noindent {\it Convex Approximation Property} (CAP).  Every holomorphic map $K\to X$ from a convex compact subset $K$ of $\C^n$ can be approximated uniformly on $K$ by holomorphic maps $\C^n\to X$.

To be clear, the notion of convexity used here is the familiar, elementary one, meaning that along with any two of its points, $K$ contains the line segment joining them.

We can also consider convex domains in $\C^k$ embedded as closed complex submanifolds of $\C^n$, $n>k$, and ask for the following version of the Cartan extension theorem.

\noindent{\it Convex Interpolation Property} (CIP).  If $T$ is a closed complex submanifold of $\C^n$ which is biholomorphic to a convex domain in some $\C^k$, then every holomorphic map $T\to X$ extends to a holomorphic map $\C^n\to X$.

CAP was introduced in \cite{FF:EOP, FF:CAP}.  CIP, introduced in \cite{Larusson5}, is a restricted version of the Oka property with interpolation (see below).  It is easily seen that CIP implies CAP (see \cite{Larusson3}), but the converse implication only comes as part of the general theory and no simple proof is known.  It turns out that CIP is unchanged if $T$ is allowed to be any holomorphically contractible submanifold of $\C^n$, or any topologically contractible submanifold of $\C^n$.

Next we will formulate a similar property for maps from an arbitrary Stein manifold $S$ to $X$.  This property implies both CAP and CIP.  We must take account of the fact that a holomorphic map from a closed complex submanifold $T$ of $S$ to $X$ need not admit a continuous extension $S\to X$.  A similar topological obstruction may appear when trying to approximate maps $K\to X$ from a compact subset $K$ in a Stein manifold $S$ by global maps $S\to X$.  By a Stein inclusion $T\hookrightarrow S$, we mean the inclusion into a Stein manifold $S$ of a closed complex submanifold $T$.

\noindent{\it Basic Oka Property} (BOP).  For every Stein inclusion $T\hookrightarrow S$ and every holomorphically convex compact subset $K$ of $S$, a continuous map $f:S\to X$ that is holomorphic on $K\cup T$ can be deformed to a holomorphic map $S\to X$.  The deformation can be kept fixed on $T$.  Also, the intermediate maps can be kept holomorphic and arbitrarily close to $f$ on $K$.

We get the implication BOP $\Rightarrow$ CAP by taking $K$ to be a convex compact subset of $S=\C^n$ and $T=\varnothing$, and we see that BOP $\Rightarrow$ CIP by choosing $K=\varnothing$ and $T\hookrightarrow S=\C^n$ as in the definition of CIP.  In these situations there are no topological obstructions to extending the map from a submanifold or approximately extending it from a convex set.

Oka properties such as BOP express a certain {\it holomorphic flexibility}.  They are opposite to Kobayashi-Eisenman-Brody hyperbolicity, which expresses {\it holomorphic rigidity}.  Recall that an $n$-dimensional complex manifold $X$ is $k$-hyperbolic in the sense of Brody for an integer $k\in \{1,\ldots, n\}$ if every holomorphic map $\C^k\to X$ is everywhere degenerate in the sense that its rank is smaller than $k$ at each point of $\C^k$.  In particular, for $k=1$, Brody hyperbolicity means that every holomorphic map $\C\to X$ is constant.  If $X$ is compact, this coincides with the notion of Kobayashi hyperbolicity \cite{Brody}. 

The theorems of Oka-Weil and of Cartan imply that $\C$, and hence $\C^n$, satisfies BOP.  Clearly, if $X$ is Brody hyperbolic, or more generally Brody volume hyperbolic, that is, Brody ($\dim X$)-hyperbolic, then $X$ does not satisfy BOP.  Moreover, no compact complex manifold $X$ of general Kodaira type satisfies BOP.  In fact, a holomorphic map $\C^{\dim X} \to X$ cannot have maximal rank at any point (easy consequence of Theorem 2 in \cite{Kobayashi-Ochiai}).

Now we introduce the parametric Oka property, which pertains to families of continuous maps parametrized by compact sets in Euclidean space.

\noindent{\it Parametric Oka Property} (POP).  Let $T\hookrightarrow S$ and $K\subset S$ be as in BOP, and let $Q\subset P$ be compact subsets of $\R^m$.  For every continuous map $f:S\times P\to X$ such that 
\begin{itemize}
\item $f(\cdot,x):S \to X$ is holomorphic for every $x\in Q$, and
\item $f(\cdot,x)$ is holomorphic on $K\cup T$ for every $x\in P$,
\end{itemize}
there is a continuous deformation $f_t:S\times P\to X$, $t\in[0,1]$, of $f=f_0$ such that
\begin{itemize}
\item the deformation is fixed on $(S\times Q)\cup (T\times P)$,
\item for every $t\in[0,1]$, $f_t(\cdot,x)$ is holomorphic on $K$ for every $x\in P$ and $f_t$ is uniformly close to $f$ on $K\times P$, and 
\item $f_1(\cdot,x):S\to X$ is holomorphic for every $x\in P$.
\end{itemize}

Ignoring the approximation condition on $K$, POP is illustrated by the following diagram.  Every lifting in the big square can be deformed through such liftings to a 
lifting in the left-hand square.  The horizontal arrows in the right-hand square
are the inclusions of the spaces of holomorphic maps into the spaces of continuous maps (these spaces carry the compact-open topology), and the vertical arrows are the restriction maps.
\[ \xymatrix{
Q \ar[r] \ar[d] &\sO(S,X) \ar[r] \ar[d] & 
\sC(S,X) \ar[d] \\ P \ar[r] \ar@{-->}[urr] \ar@{..>}[ur] 
& \sO(T,X) \ar[r] & \sC(T,X)
} \]

Clearly, BOP is the special case of POP obtained by taking $P$ to be a singleton and $Q$ empty.

The following classical result, due to H.\ Grauert, is the basis of the {\it Oka-Grauert principle}. 

\begin{theorem}[Grauert \cite{Grauert3}]
Every complex homogeneous manifold satisfies {\rm POP} for all pairs of finite polyhedra $Q\subset P$.  The analogous result holds for sections of holomorphic 
$G$-bundles, for any complex Lie group $G$, over a Stein space.
\end{theorem}

Since an isomorphism between $G$-bundles is represented by a section of an associated $G$-bundle, we get the following corollary.

\begin{corollary}[The Oka-Grauert principle for vector bundles]
The holomorphic and topological classifications of such bundles over Stein spaces coincide.  This holds in particular for complex vector bundles.
\end{corollary}

K.\ Oka, a pioneer of several complex variables after whom the Oka principle is named, proved this result in 1939 for line bundles over domains of holomorphy \cite{Oka1939}.  He showed that on a domain of holomorphy, the second Cousin problem is solvable by holomorphic functions if it is solvable by continuous functions.

Among the Oka properties introduced above, CAP and CIP are ostensibly the weakest, and hence the easiest to verify for concrete examples, while POP is ostensibly the strongest. It is surprising and far from obvious that these properties are in fact mutually equivalent.

\begin{theorem}
\label{equiv}
For every complex manifold $X$, each Oka property implies the others:
\[	\rm{CAP} \Leftrightarrow \rm{CIP} \Leftrightarrow 
	{\rm BOP} \Leftrightarrow {\rm POP}.  \] 
These equivalences also hold if the inclusion $T\hookrightarrow S$ in {\rm BOP} and {\rm POP} is allowed to be the inclusion in a reduced Stein space of a closed analytic subvariety.
\end{theorem}

\begin{definition}
\label{Okamanifold}
A complex manifold satisfying the equivalent properties CAP, CIP, BOP, and POP is called an {\it Oka manifold}.
\end{definition}

The question whether it is possible to characterize BOP and POP by a Runge approximation property for maps $\C^n\to X$ was raised by M.\ Gromov in \cite{Gromov:OP}. The implication CAP $\Rightarrow$ BOP was established in \cite{FF:EOP, FF:CAP} (see \cite{FF:Kohn} for singular Stein source spaces).  In the same papers it was shown that a certain parametric version of CAP implies POP.  The final implication BOP $\Rightarrow$ POP was established in \cite{FF:OkaManifolds}, where the class of Oka manifolds was first formally introduced.

Our definition of BOP corresponds to what has been called in the literature the {\em basic Oka property with approximation and interpolation} (BOPAI).  By removing the interpolation condition in BOP we get the {\em basic Oka property with approximation} (BOPA); similarly, by removing the approximation condition we get  
the {\em basic Oka property with interpolation} (BOPI).  Analogous statements hold for the parametric Oka properties.  It was observed in \cite{Larusson3} that BOPI $\Rightarrow$ BOPA and POPI $\Rightarrow$ POPA (not only for manifolds, but also for holomorphic maps; see \S\ref{maps} below).  This was the first nontrivial implication to be noted between variants of the Oka property.  We now know that all these variants are equivalent.

\smallskip\noindent 
{\bf Properties of Oka manifolds.} 
It is clear from CAP that an Oka manifold $X$ is {\it strongly dominable} by $\C^n$, $n=\dim X$, in the sense that for every $x\in X$, there is a holomorphic map $f:\C^n\to X$ such that $f(0)=x$ and $f$ is a local biholomorphism at $0$.  Hence, the Kobayashi pseudometric of $X$ vanishes identically.  Also, if a plurisubharmonic function on a connected Oka manifold is bounded above, then it is constant.  It follows that every $\R$-complete holomorphic vector field on a Stein Oka manifold is $\C$-complete \cite{FF:Actions}.

A deeper property of Oka manifolds is jet transversality for holomorphic maps from Stein manifolds.  If $S$ is a Stein manifold, $X$ is an Oka manifold, $k\geq 0$, and $B$ is a closed complex submanifold of the complex manifold of $k$-jets of holomorphic maps $S\to X$, then the $k$-jet of a generic holomorphic map $S\to X$ is transverse to $B$.  (This is a simplified version of Theorem 4.2 in \cite{FF:Flexibility}.)  Here, a property that holds for all maps in a countable intersection of open dense subsets of $\mathscr O(S,X)$ with the compact-open topology is said to be generic.  Without the Oka assumption on $X$, S.~Kaliman and M.~Zaidenberg have shown that if $f:S\to X$ is holomorphic and $K\subset S$ is compact, then there is a holomorphic map $K\to X$, uniformly approximating $f|K$ as closely as desired, whose $k$-jet is transverse to $B$ \cite{Kaliman-Zaidenberg}.

The jet transversality theorem can be used to show that if $S$ is a Stein manifold and $X$ is an Oka manifold, then a generic holomorphic map $S\to X$ is an immersion when $\dim X\geq 2\dim S$, and an injective immersion when $\dim X\geq 2\dim S+1$ (\cite{FF:Flexibility}, Corollary 1.5).

Applying POP with the parameter space pairs $\varnothing\hookrightarrow\ast$, $\{0,1\}\hookrightarrow [0,1]$, $\ast\hookrightarrow S^k$ and $S^k\hookrightarrow B^{k+1}$ for $k\geq 1$, where $B^{k+1}$ denotes the $(k+1)$-dimensional ball and $S^k=\partial B^{k+1}$ is the $k$-dimensional sphere, we obtain the following result.

\begin{corollary}
For a Stein manifold $S$ and an Oka manifold $X$, the inclusion 
\[\sO(S,X)\hookrightarrow\sC(S,X)\]
is a weak homotopy equivalence with respect to the compact-open topology, that is, it induces a bijection of path components and isomorphisms of all homotopy groups.
\end{corollary} 

The following functorial property of Oka manifolds was proved in \cite{FF:CAP, FF:OkaManifolds}.

\begin{theorem} 
\label{up-and-down}
If $E$ and $B$ are complex manifolds and $E\to B$ is a holomorphic fibre bundle
whose fibre is an Oka manifold, then $B$ is an Oka manifold if and only if $E$ is an Oka manifold. 
\end{theorem}

We record separately the important special case when $E\to B$ is a covering map.

\begin{corollary}
\label{covers-quotients}
If $E\to B$ is a holomorphic covering map, then $B$ is an Oka manifold if and only if $E$ is an Oka manifold.  
\end{corollary}

\noindent{\bf Examples of Oka manifolds.}  We end this section by listing all the explicit examples of Oka manifolds that are currently known.

\begin{itemize}
\item  The Riemann surfaces that are Oka are $\C$, $\C^*$, the Riemann sphere, and all tori.  In other words, a Riemann surface is Oka if and only if it is not hyperbolic, that is, not covered by the disc.
\item  Complex affine spaces $\C^n$, complex projective spaces $\P^n$ for all $n$; Grassmannians.
\item  More generally, complex Lie groups and their homogeneous spaces.
\item  $\C^n\backslash A$, where $A$ is an algebraic subvariety or a tame analytic subvariety of complex codimension at least $2$.
\item  $\P^n\backslash A$, where $A$ is an algebraic subvariety 
of codimension at least $2$.
\item  Hirzebruch surfaces (holomorphic $\P^1$-bundles over $\P^1$).
\item  Hopf manifolds (compact manifolds with universal covering space $\C^n\backslash\{0\}$, $n\geq 2$).
\item  Algebraic manifolds that are Zariski locally affine (every point has a Zariski neighbourhood isomorphic to $\C^n$); modifications of such obtained by blowing up points or removing subvarieties of codimension at least $2$.
\item  $\C^n$ blown up at all points of a tame discrete sequence.
\item  Complex tori of dimension at least $2$ with finitely many points removed, or blown up at finitely many points.
\end{itemize}

Most of these examples are elliptic complex manifolds (see \S\ref{elliptic} below); none are known not to be elliptic.  All are at least weakly subelliptic.  

The most important method to produce new Oka manifolds from old, and thus expand the collection of examples, is provided by the result (a generalization of Theorem \ref{up-and-down}) that if $E$ and $B$ are complex manifolds and $E\to B$ is a surjective Oka map (see \S\ref{maps} below), then $E$ is Oka if and only if $B$ is Oka.


\section{The proof that CAP implies POP}
\label{proof}

\noindent
The proof of the main implication CAP $\Rightarrow$ POP in Theorem \ref{equiv} is broken down into two steps, CAP $\Rightarrow$ BOP and BOP $\Rightarrow$ POP.  The first implication was proved in \cite{FF:EOP, FF:CAP}; we shall outline here the main steps. The second implication was established in \cite{FF:OkaManifolds} using a recent theorem from \cite{Forstneric-Wold2} on the existence of Stein neighbourhood bases of certain compact sets in $\C^n\times\C^m$ that are fibred over the real space $\R^n \subset \C^n$ such that the fibres are Stein compacts.   

The scheme of the proof of CAP $\Rightarrow$ BOP is essentially that used in the proof of Cartan's Theorems A and B or, more specifically, in the proof of the Oka-Grauert principle given by Henkin and Leiterer \cite{Henkin-Leiterer:Oka}.  Assume that $\pi:Z\to S$ is a holomorphic fibre bundle over a Stein manifold $S$, whose fibre $X$ satisfies CAP.  We exhaust $S$ by an increasing sequence of compact strongly pseudoconvex domains $A_k$ that are $\sO(S)$-convex.  The domain $A_{k+1}=A_k\cup B_k$ is obtained by attaching to $A_k$ either a convex bump $B_k$ or a special handle whose core is a totally real disc attached to $\partial A_k$ along a complex tangential sphere.  (The sets $A_k$ are obtained using a strictly plurisubharmonic exhaustion $\rho:S\to \R$.  We may choose $A_0$ to be a small neighbourhood of the compact set $K=\widehat K_{\sO(S)}\subset S$ on which the initial continuous section $f_0:S\to Z$ is holomorphic.)  Each $B_k$ can be chosen small enough that our fibre bundle $Z\to S$ is trivial over a neighbourhood of it.  A holomorphic section $f:S\to Z$, homotopic to the initial section $f_0$, is found as a locally uniform limit $f=\lim\limits_{k\to\infty} f_k$ of a sequence of continuous sections $f_k:S\to Z$ such that $f_k$ is holomorphic over a neighbourhood of $A_k$.  In the induction step, we find $f_{k+1}$ that approximates $f_k$ uniformly on a neighbourhood of $A_k$ and is homotopic to it.  We treat separately the extension across a convex bump (the {\it noncritial case}) and the crossing of a critical level of an exhaustion function on $X$ (the {\it critical case}).

In the {\it noncritical case} (which is the more difficult of the two, and is the only one where the CAP of the fibre $X$ is used) the problem is divided in two substeps.  First we approximate $f_k$ on a neighbourhood of the set $C_k=A_k\cap B_k$ by a holomorphic section $g_k$ defined on a neighbourhood of $B_k$.  Since the bundle is trivial there, this is a Runge approximation problem for maps to the fibre $X$, and here the CAP of the fibre is invoked.  We then glue $f_k$ and $g_k$ into a section $f_{k+1}$ which is holomorphic over a neighbourhood of $A_{k+1}=A_k\cup B_k$.

Let us describe this procedure a bit more carefully.  In the classical case when $f_k$ and $g_k$ are maps to a complex Lie group $G$ (the fibre of a principal holomorphic bundle $Z\to S$) one has $f_k=g_k\cdot \gamma_k$ on $C_k$, where $\gamma_k = g_k^{-1}\cdot f_k:C_k\to G$ is a holomorphic map with values in $G$ that is close to the constant map $x\mapsto 1\in G$ (the identity element of $G$).  By the classical Cartan lemma we can split $\gamma_k$ into a product $\gamma_k=\beta_k\cdot\alpha_k^{-1}$ with holomorphic maps $\alpha_k:A_k\to G$, $\beta_k:B_k\to G$ that are close to the constant map $x\mapsto 1$, and we take $f_k\cdot \alpha_k=g_k\cdot \beta_k$ as the next map $f_{k+1}$.

In the general case we work with {\it thick sections}, also called {\it local holomorphic sprays}.  These are families of holomorphic sections, depending holomorphically on a parameter in a neighbourhood of the origin in a Euclidean space, which are submersive with respect to the parameter.  Given a thick section $f_k$ over $A_k$, we approximate $f_k$ uniformly over a neighbourhood of $C_k$ by a thick section $g_k$ over $B_k$ as described above, using the CAP of the fibre $X$.  We then find a fibre-preserving biholomorphic transition map $\gamma_k$, close to the identity map and satisfying 
\[
	f_k=g_k\circ \gamma_k \quad \hbox{near } C_k.
\]
Next we split $\gamma_k$ into a composition of the form
\[ 
	\gamma_k = \beta_k \circ\alpha_k^{-1},
\]
where $\alpha_k$ and $\beta_k$ are holomorphic maps over $A_k$ and $B_k$, respectively.  It follows as before that 
\[	
	f_k\circ \alpha_k = g_k\circ \beta_k \quad \hbox{over } C_k,
\]
and hence these two thick sections amalgamate into a thick section $f_{k+1}$ that is holomorphic over a neighbourhood of $A_{k+1}$.  The induction may proceed.

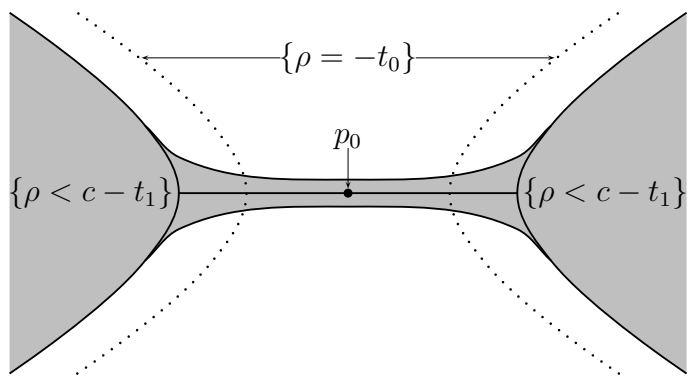
\begin{figure}[ht]
\psset{unit=0.6cm, xunit=1.5, linewidth=0.7pt}
\begin{pspicture}(-6.5,-4.3)(6.5,4.5)


\pscustom[fillstyle=solid,fillcolor=lightgray,linestyle=none]  
{\pscurve(-3,-1.5)(-2.5,-0.8)(0,-0.3)(2.5,-0.8)(3,-1.5) 
\psline[linestyle=dashed,linewidth=0.2pt](3,-1.5)(3,1.5)
\pscurve[liftpen=1](3,1.5)(2.5,0.8)(0,0.3)(-2.5,0.8)(-3,1.5)
\psline[linestyle=dashed,linewidth=0.2pt](-3,1.5)(-3,-1.5)}

\pscustom[fillstyle=solid,fillcolor=lightgray]
{\pscurve[liftpen=1](5,4)(3,1.5)(2.5,0)(3,-1.5)(5,-4)}  

\pscustom[fillstyle=solid,fillcolor=lightgray]
{\pscurve[liftpen=1](-5,4)(-3,1.5)(-2.5,0)(-3,-1.5)(-5,-4)}  


\psline(-2.5,0)(2.5,0) 
\psecurve(5,4)(3,1.5)(2.5,0.8)(0,0.3)(-2.5,0.8)(-3,1.5)(-5,4) 
\psecurve(5,-4)(3,-1.5)(2.5,-0.8)(0,-0.3)(-2.5,-0.8)(-3,-1.5)(-5,-4)  


\pscurve[linestyle=dotted,linewidth=1pt](4,4)(2,1.5)(1.5,0)(2,-1.5)(4,-4)             \pscurve[linestyle=dotted,linewidth=1pt](-4,4)(-2,1.5)(-1.5,0)(-2,-1.5)(-4,-4)
         

\rput(0,3){$\{\rho=-t_0\}$}
\psline[linewidth=0.2pt]{->}(1,3)(3.05,3)
\psline[linewidth=0.2pt]{->}(-1,3)(-3.05,3)



\psdot(0,0)
\rput(0,1.2){$p_0$}
\psline[linewidth=0.2pt]{->}(0,1)(0,0.1)

\rput(3.8,0){$\{\rho<c-t_1\}$}
\rput(-3.8,0){$\{\rho<c-t_1\}$}

\end{pspicture}
\caption{Passing a critical level of $\rho$.}
\label{Fig:tau}
\end{figure}

The {\it critical case} is treated by approximately extending the holomorphic section across the stable manifold at the critical point $p_0$, thereby reducing the 
problem to the noncritical case for a different strictly plurisubharmonic function that is used only to cross the critical level $\rho=\rho(p_0)$; we then revert back to the original function $\rho$. The relevant geometry near the critical level is illustrated in Fig.\ \ref{Fig:tau}, where we assume that $\rho(p_0)=0$.

In the case of a stratified fibre bundle over a possibly singular Stein space $S$, and when interpolating a section on a subvariety $T$ of $S$, essentially the same proof can be accomplished by induction on the strata (see \cite{FF:Kohn}).

Every step of the proof can also be carried out in the parametric case, and this shows that a parametric version of CAP (called PCAP) for the fibres of $Z\to S$ implies the parametric Oka principle for sections $S\to Z$. The proof of Theorem \ref{equiv} is completed by showing that POP $\Rightarrow$ PCAP (see \cite{FF:OkaManifolds}).


\section{Elliptic and subelliptic manifolds}
\label{elliptic}

\noindent
Grauert established the Oka principle for sections of a holomorphic fibre bundle over a Stein manifold if the fibres are homogeneous spaces with respect to a complex Lie group $G$ and the bundle has a trivialization whose transition functions take values in $G$.  The key innovation in Gromov's paper \cite{Gromov:OP} is the generalization from homogeneous manifolds to elliptic manifolds and the realization that the Oka principle holds for sections of a holomorphic fibre bundle with elliptic fibres over a Stein manifold regardless of how the bundle may be trivialized.

\begin{definition}[Gromov \cite{Gromov:OP}]  A {\it spray} on a complex manifold $X$ is a holomorphic map $s:E\to X$ defined on the total space of a holomorphic vector bundle $E$ over $X$ such that $s(0_x)=x$ for all $x\in X$.  The spray $s$ is said to be {\it dominating at} $x\in X$ if $s|E_x\to X$ is a submersion at $0_x$, and $s$ is said to be {\it dominating} if it is dominating at every point of $X$.  Finally, $X$ is said to be {\it elliptic} if it admits a dominating spray.
\end{definition}

The following ostensibly more general notion, defined in \cite{FF:Subelliptic}, is sometimes easier to verify.

\begin{definition}  A complex manifold $X$ is said to be {\it subelliptic} if it admits finitely many sprays $s_j:E_j\to X$ that together dominate at every $x\in X$, that is, the images $(ds_j)_{0_{j,x}}(E_{j,x})$ together span the tangent space $T_x X$ (here, we identify the vector space $E_{j,x}$ with its tangent space at the zero element $0_{j,x}$).
\end{definition}

A further weakening of ellipticity, called {\it weak subellipticity}, requires the existence of countably many sprays that together dominate at each point.

\begin{figure}[ht]
\psset{unit=0.6cm,linewidth=1pt}

\begin{pspicture}(0,0)(6,6.5)

\pscustom[fillstyle=solid,fillcolor=lightgray]
{
\pscurve(1,1)(4,2)(6.5,2)
\pscurve(6.5,2)(6.5,3.8)(7.3,5)
\pscurve(7.3,5)(5,5)(2,4)
\pscurve(2,4)(1.8,2.6)(1,1)
}

\psdot[dotsize=3pt](4,3.5)
\psdot[dotsize=3pt](3.5,5)
\psline(4,3.5)(3,6.5)
\psline[linestyle=dotted](4,3.5)(4.5,2)
\psline(4.5,2)(5,0.5)
\psarc[linewidth=0.4pt,arrows=<-](4,3.5){1.6}{30}{110}
\psdot[dotsize=3pt](5.38,4.3)

\rput(4.5,5.4){$s$}
\rput(2.5,2.5){$X$}
\rput(3.7,3.2){$x$}
\rput(1,5){$E_x$}
\psline[linewidth=0.21pt]{->}(1.3,5.2)(3.1,5.8)

\end{pspicture}
\caption{The fibre $E_x$ of a spray on the manifold $X$.}
\end{figure}
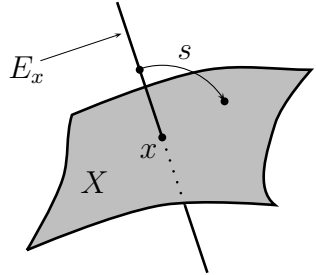

\noindent
{\bf Examples.}  Here are the most familiar examples of elliptic and subelliptic manifolds.  The first three appeared in \cite{Gromov:OP}.

\noindent 
$\bullet$ \ A homogeneous space of a complex Lie group $G$, that is, a complex manifold $X$ on which $G$ acts holomorphically and transitively, is elliptic.  The map $s:X\times\mathfrak g\to X$, $(x,v)\mapsto \exp(v)\cdot x$, where $\mathfrak g$ is the Lie algebra of $G$, is a dominating spray defined on a trivial vector bundle over $X$.  In particular, a complex Lie group is elliptic.

\noindent 
$\bullet$ \ More generally, if $X$ admits finitely many $\C$-complete holomorphic vector fields $v_1,\ldots,v_k$ that span $T_x X$ at every point $x\in X$, then the map $s:X\times \C^k\to X$, 
\[ s(x,t_1,\ldots,t_k)= \phi_{t_1}^1\circ\cdots\circ \phi_{t_k}^k (x), \]
where $\phi^j_t$ is the flow of $v_j$, is a dominating spray on $X$.

\noindent 
$\bullet$ \ A spray of this type exists on $X=\C^n\backslash A$ where $A$ is an algebraic subvariety of co\-dimension at least $2$.  We take the fields $v_1,\dots, v_k$ to be of the form $v(z)=f(\pi(z))b$, where $b\in\C^n\setminus\{0\}$, $\pi:\C^n\to\C^{n-1}$ is a linear projection with $\pi(b)=0$ such that $\pi|A$ is proper, and $f:\C^{n-1}\to\C$ is a polynomial that vanishes on the subvariety $\pi(A)$.  The flow of $v$, given by $\phi_t(z)=z+tf(\pi(z))b$, fixes $A$ and thus restricts to a complete flow on $\C^n\setminus A$.  

\noindent
$\bullet$ \ If $A$ is a subvariety of $\P^n$ of codimension at least $2$, then $\P^n\setminus A$ is subelliptic (\cite{FF:Subelliptic}, Proposition 1.2).  It is unknown whether these manifolds are elliptic.

\smallskip
The product of elliptic manifolds is elliptic.  The only other known method to obtain new elliptic manifolds from old is provided by the simple observation that if $X\to Y$ is an unbranched holomorphic covering map and $Y$ is elliptic, then $X$ is elliptic.  The same is true of subellipticity.  By contrast, the Oka property is known to pass up and down along much more general maps.

The central results in Gromov's 1989 paper are the following.

\begin{theorem}[Gromov \cite{Gromov:OP}]  The parametric Oka property holds, in order of increasing generality, for:
\begin{enumerate}
\item[(i)]  maps from a Stein manifold to an elliptic manifold (so an elliptic manifold is Oka),
\item[(ii)]  sections of holomorphic fibre bundles with elliptic fibres over a Stein manifold, and
\item[(iii)]  sections of elliptic submersions over a Stein manifold.
\end{enumerate}
\end{theorem}

The definition of the parametric Oka property for sections is analogous to that given in \S 3 for maps.  A holomorphic submersion $h:Z\to X$ between reduced complex spaces is said to be {\it elliptic} if every point of $X$ has an open neighbourhood $U$ such that the restriction $h^{-1}(U)\to U$ admits a dominating fibre spray.  A {\it fibre spray} for a holomorphic submersion $h:Z\to X$ is a holomorphic vector bundle $\pi:E\to Z$ with a holomorphic map $s:E\to Z$ such that $h\circ\pi=h\circ s$ and $s(0_z)=z$ for each $z\in Z$.  We say that the fibre spray is {\it dominating at} $z\in Z$ if the derivative of $s$ at $0_z$ maps the subspace $E_z$ of $T_{0_z}E$ onto the vertical subspace $\ker dh_z$ of $T_z Z$.  We say that the fibre spray is {\it dominating} if it is dominating at every point.  A {\it subelliptic submersion} is defined similarly.  Clearly, a holomorphic fibre bundle with (sub)elliptic fibres is a (sub)elliptic submersion.

Detailed proofs of Gromov's theorem were given by J.\ Prezelj and the first-named author in \cite{Forstneric-Prezelj2000, Forstneric-Prezelj2001, Forstneric-Prezelj2002}.  Gromov's theorem was generalized to subelliptic submersions in \cite{FF:Subelliptic}.

\smallskip \noindent
{\bf Example.} Let $\pi:E\to S$ be a holomorphic vector bundle with fibre $\C^k$ over a Stein space $S$.  Let $\widehat E\to S$ be the associated bundle with fibre $\P^k$.  Let $\Sigma\subset E$ be a subvariety with fibres of codimension at least $2$ over each point of $S$, such that the closure of $\Sigma$ in $\widehat E$ does not contain the hyperplane at infinity over any point of $S$.  Then $E\setminus\Sigma\to S$ is an elliptic submersion.  Thus, sections $S\to E$ of $\pi$ avoiding $\Sigma$ satisfy the parametric Oka property.

\smallskip
Gromov's theorem has been generalized to stratified maps as follows.

\begin{theorem}[\cite{FF:Kohn}]
The parametric Oka property holds for sections of stratified holomorphic fibre bundles with Oka fibres and for sections of stratified subelliptic submersions over reduced Stein spaces.
\end{theorem}

A {\it stratification} of a complex space $X$ is a descending chain $X=X_0\supset X_1 \supset\cdots\supset X_m=\varnothing$ of closed subvarieties such that for $k=1,\dots,m$, each connected component $S$ of $X_{k-1}\setminus X_k$ (that is, each {\it stratum}) is smooth and $\overline{S}\setminus S\subset X_k$.  A holomorphic submersion $\pi:Z\to X$ is a {\it stratified holomorphic fibre bundle} if $X$ admits a stratification such that for every stratum $S$, the restriction $\pi^{-1}(S)\to S$ is a holomorphic fibre bundle over $S$.  The fibres over different strata may be different.  Also, $\pi$ is a {\it stratified (sub)elliptic submersion} if $X$ admits a stratification such that for every stratum $S$, the restriction $\pi^{-1}(S)\to S$ is (sub)elliptic.

The strongest available version of Gromov's theorem is Theorem \ref{strongest-theorem} below.

\smallskip\noindent 
{\bf Applications.}  The following are some of the applications of Gromov's Oka principle.

\noindent
$\bullet$ \ The proof of Forster's conjecture was mentioned in \S 2.

\noindent
$\bullet$ \ The following homotopy principle for holomorphic immersions is due to 
Eliashberg and Gromov (see \cite{Eliashberg-Gromov1} and \cite{Gromov:PDR}, \S 2.1.5).  If $X$ is a Stein manifold whose complex cotangent bundle $T^*X$ is generated by $(1,0)$-forms $\theta_1,\dots,\theta_q$ where $q>\dim X$, then the $q$-tuple $(\theta_1,\dots,\theta_q)$ can be continuously deformed through $q$-tuples generating $T^*X$ to a $q$-tuple $(df_1,\dots,df_q)$, where $(f_1,\dots,f_q):X\to\C^q$ is a holomorphic immersion.  Consequently, every $n$-dimensional Stein manifold admits a holomorphic immersion into $\C^{[3n/2]}$.  D.\ Kolari\v c has given a complete exposition of Eliashberg and Gromov's result, as well as a $1$-parametric version of it \cite{Kolaric}.

\noindent
$\bullet$ \ There is also a homotopy principle for holomorphic submersions \cite{FF:Noncritical}.  While not actually an application of Gromov's Oka principle, it relies on methods that have emerged in the course of its development.  Let $X$ be a Stein manifold and $1\leq q<\dim X$.  By a $q$-{\it coframe} on $X$ we mean a $q$-tuple of continuous $(1,0)$-forms on $X$ that are linearly independent at each point.  Every $q$-coframe on $X$ can be continuously deformed through $q$-coframes to a $q$-coframe $(df_1,\dots,df_q)$, where $(f_1,\dots,f_q):X\to\C^q$ is a holomorphic submersion.  Consequently, every $n$-dimensional Stein manifold admits a holomorphic submersion into $\C^k$, where $k=[(n+1)/2]$.  This value of $k$ is sharp for every $n$.  Also, an $n$-dimensional parallelizable Stein manifold, $n\geq 2$, has a submersion into $\C^{n-1}$.  Whether $\C^{n-1}$ can be replaced by $\C^n$ is an open problem of great interest.

\noindent
$\bullet$ \ The following theorem is a special case of the results of \cite{FF:Intersections} on the elimination of intersections.  Let $A$ be an algebraic subvariety of $\C^d$ of codimension at least $2$.  Let $X$ be a Stein manifold, $f:X\to \C^d$ be holomorphic, and $Y\subset X$ be a (possibly empty) union of connected components of $f^{-1}(A)$.  If there is a continuous map $\tilde f:X\to\C^d$ which equals $f$ in a neighbourhood of $Y$ and satisfies $\tilde f^{-1}(A)=Y$, then for each $r\in\mathbb N$, there is a holomorphic map $g:X\to\C^d$ such that $g^{-1}(A)=Y$ and $g-f$ vanishes to order $r$ along $Y$.  Such a map $g$ always exists if $\dim X<2(d-\dim A)$, or if $X$ is contractible and $\dim Y\leq 2(d-\dim A-1)$.

\noindent
$\bullet$ \ The final application we shall describe is the recent solution of the Gromov-Vaserstein problem by B.\ Ivarsson and F.\ Kutzschebauch.  The problem of factoring a matrix in $\mathrm{SL}(m,R)$, where $R$ is a commutative ring, into a product of elementary matrices has been studied for many rings.  For the ring of complex numbers, it is of course a basic result of linear algebra that this is always possible.  For more sophisticated rings, such as rings of polynomials or continuous functions, the problem has been considered by a number of mathematicians, including A.~A.~Suslin, W.\ Thurston, and L.\ Vaserstein.  A holomorphic version of the problem was posed by Gromov in \cite{Gromov:OP}, 3.5.G, and solved as follows.

\begin{theorem}[Ivarsson and Kutzschebauch \cite{Ivarsson-Kutzschebauch}]
\label{IK}
Let $f:S\to \mathrm{SL}(m,\C)$ be a null\-homotopic holomorphic map from a Stein manifold $S$.  Then there are holomorphic maps $G_1,\dots, G_k:S\to \mathbb{C}^{m(m-1)/2}$ such that
\[		f = \left[\begin{array}{cc} 1  &  0 \cr G_1 & 1 \cr \end{array} \right]  
		\left[\begin{array}{cc} 1 & G_2 \cr 0 & 1 \cr \end{array} \right]
		\left[\begin{array}{cc} 1 & 0 \cr G_3 & 1 \cr \end{array} \right] \cdots.  \]
\end{theorem}

The algebraic case was considered much earlier.  In 1966, P.\ M.\ Cohn showed that the matrix
\[	\left[\begin{array}{cc} 1-zw &  z^2\\ -w^2 & 1+zw \end{array}\right]
	 \in \mathrm{SL}(2,\mathbb{C}[z,w])  \]
does not factor as a finite product of unipotent matrices.  Yet, by Theorem \ref{IK}, it is a product of unipotent matrices over the larger ring $\sO(\C^2)$.  In 1977, Suslin proved that for $m\geq 3$ (and any $n$), every matrix in $\mathrm{SL}(m,\C{[z_1,\dots,z_n]})$ is a product of unipotent matrices.

\smallskip\noindent
{\it Idea of proof of Theorem \ref{IK}:}  Define $\Psi_k:(\C^{m(m-1)/2})^k\to \mathrm{SL}(m,\C)$ by
\[  \Psi_k(g_1,\dots,g_k) = \left[\begin{array}{cc} 1  &  0 \cr g_1 & 1 \cr \end{array} \right] 
         \left[\begin{array}{cc} 1 & g_2 \cr 0 & 1 \cr \end{array} \right]    	
	 \left[\begin{array}{cc} 1 & 0\cr g_3 & 1 \cr \end{array} \right]\cdots.  \]
We want to find a holomorphic map $G:S\to (\C^{m(m-1)/2})^k$ such that the following diagram commutes:
\[	\xymatrix{ & (\C^{m(m-1)/2})^k \ar[d]^{\Psi_{k}} \\ 
	 S \ar[r]_{\!\!\!\!\!\!\!\!\!f} \ar[ur]^{G} & \mathrm{SL}(m,\C)}  \]
A theorem of Vaserstein gives a continuous lifting of $f$.  Unfortunately, the map $\Psi_k$ is not a submersion and its fibres are hard to analyze.  Still, the continuous lifting may be deformed to a holomorphic lifting by applying the Oka principle to certain related maps that are stratified elliptic submersions.  The version of Gromov's Oka principle proved in \cite{FF:Kohn} is required.  \hfill\qed

\smallskip
Ivarsson and Kutzschebauch point out that as a consequence of their theorem, the inclusion of the ring of holomorphic functions on a contractible Stein manifold into the ring of continuous functions does not induce an isomorphism of $K_1$-groups, whereas by Grauert's Oka principle it does induce an isomorphism of $K_0$-groups.  Thus, here, Gromov's Oka principle reveals a limitation of a more general (and vague) Oka principle. 

\smallskip\noindent
{\bf Variants of Gromov's Oka principle.} \ $\bullet$ \ It is possible to push the Oka principle beyond the realm of Stein spaces.  J.\ Prezelj has proved a version of Gromov's Oka principle for sections of certain holomorphic submersions over reduced 1-convex spaces \cite{Prezelj2010}.  Recall that a complex space is 1-convex if it contains a largest compact subvariety that can be blown down so as to make a Stein space.

\noindent
$\bullet$ \  As mentioned above, a holomorphic function on a subvariety $A$ of a Stein manifold $X$ extends to all of $X$. In the absence of topological obstructions, the analogous result holds for holomorphic maps to an Oka manifold $Y$, but fails for arbitrary $Y$ (for example, let $X=\C$ with a two-point subvariety and $Y$ be a disc).  M.\ Slapar and the first-named author have shown that, unless $\dim X=2$, a version of Cartan's theorem can be obtained for arbitrary $Y$ if we are permitted to deform the Stein structure on $X$ away from $A$ \cite{Forstneric-Slapar2007a, Forstneric-Slapar2007b}.  When $\dim X=2$, this still holds if we are allowed to deform not only the complex structure but also the smooth structure of $X$ away from $A$.


\section{From manifolds to maps}
\label{maps}

\noindent
Once a new property has been defined for objects, we should extend it to arrows (or at least try).  For the Oka property, this was first done in \cite{Larusson2}.  The following definition is stronger than the one in \cite{Larusson2} in that it includes approximation, but in fact the two definitions can be shown to be equivalent.

\begin{definition}
\label{POP-for-maps}
A holomorphic map $\pi:E\to B$ between complex manifolds has the {\it parametric Oka property} if whenever
\begin{itemize}
\item  $S$ is a Stein manifold,
\item  $T$ is a closed complex submanifold of $S$,
\item  $K$ is a holomorphically convex compact subset of $S$,
\item  $P$ is a finite polyhedron with a subpolyhedron $Q$,
\item  $f:S\times P\to B$ is continuous and holomorphic along $S$,
\item  $g_0:S\times P\to E$ is a continuous lifting of $f$ by $\pi$ (so $\pi\circ g_0=f$) such that $g_0(\cdot,x)$ is holomorphic on $S$ for all $x\in Q$ and holomorphic on $K\cup T$ for all $x\in P$,
\end{itemize}
there is a continuous deformation $g_t:S\times P\to E$ of $g_0$ such that for all $t\in[0,1]$,
\begin{itemize}
\item  $\pi\circ g_t=f$,
\item  $g_t=g_0$ on $S\times Q$ and on $T\times P$,
\item  $g_t|K\times P$ is holomorphic along $K$ and uniformly close to $g_0|K\times P$, and
\item  $g_1:S\times P\to E$ is holomorphic along $S$.
\end{itemize}
If the above holds when $P$ is a singleton and $Q$ is empty, then $\pi$ is said to have the {\it basic Oka property}.
\end{definition}

\begin{theorem}[\cite{FF:OkaMaps}, Theorem 1.1]
\label{basic=parametric}
For holomorphic submersions, the basic Oka property implies the parametric Oka property.
\end{theorem}

As the definition is somewhat complicated, it may help to consider the very simplest case when $P$ is a singleton and $Q$, $T$, and $K$ are empty.  Then the Oka property of $\pi$ simply says that if $S$ is Stein and $f:S\to B$ is holomorphic, then every continuous lifting $S\to E$ of $f$ by $\pi$ can be continuously deformed through such liftings to a holomorphic lifting.  In particular, if $B$ is Stein, every continuous section of $\pi$ can be continuously deformed to a holomorphic section.

The class of maps with the parametric Oka property is not closed under composition, as shown by the simple example $\mathbb D=\{z\in\C: \lvert z\rvert<1\}\hookrightarrow\C\to\ast$.  The following notion, motivated by homotopy-theoretic considerations (see \S 7), is more natural.

\begin{definition}
A holomorphic map $\pi:E\to B$ between complex manifolds is an {\it Oka map} if it is a topological fibration and satisfies the parametric Oka property.
\end{definition}

\noindent
{\bf Remarks.} \ $\bullet$ \ By a topological fibration we mean a Serre fibration or a Hurewicz fibration: the two notions are equivalent for continuous maps between locally finite simplicial complexes (\cite{Arnold}, Theorem 4.4), such as complex manifolds or complex spaces.

\noindent
$\bullet$ \ A complex manifold $X$ is Oka if and only if the constant map $X\to\ast$ is Oka.

\noindent
$\bullet$ \ The class of Oka maps is closed under composition.  The pullback of an Oka map by an arbitrary holomorphic map is Oka.  A retract of an Oka map is Oka.

\noindent
$\bullet$ \ An Oka map $\pi:E\to B$ is a holomorphic submersion.  Namely, since $\pi$ is a topological fibration, the image of $\pi$ is a union of path components of $B$ and $\pi$ has a continuous section on a neighbourhood, say biholomorphic to a ball, of each point $q$ of the image, taking $q$ to any point $p$ in its preimage.  The Oka property of $\pi$ implies that such a section can be deformed to a holomorphic section still mapping $q$ to $p$.  Thus, by the holomorphic rank theorem, $\pi$ is a holomorphic submersion.

\noindent
$\bullet$ \ The fibres of an Oka map are Oka manifolds.  However, a holomorphic submersion with Oka fibres need not be an Oka map, even if it is smoothly locally trivial (see Example \ref{eremenko} below).

\noindent
$\bullet$ \ In Definition \ref{POP-for-maps}, $Q\hookrightarrow P$ may be taken to be any cofibration between cofibrant topological spaces (\cite{Larusson2}, \S 16), such as arbitrary CW-complexes, not necessarily compact.  Here, the notion of cofibrancy for topological spaces and continuous maps is the stronger one that goes with Serre fibrations rather than Hurewicz fibrations.

\noindent
$\bullet$ \ The theory of Oka maps can to a large extent be formulated for holomorphic submersions between reduced complex spaces.  Also, for some purposes, $P$ and $Q$ in Definition \ref{POP-for-maps} may be taken to be arbitrary compact subsets of $\R^m$ or even arbitrary compact Hausdorff spaces, and $S$ and $T$ may be taken to be reduced Stein spaces (see \cite{FF:OkaMaps, FF:Invariance}).  

\smallskip

The following theorem is pleasantly analogous to a basic property of Hurewicz fibrations.

\begin{theorem}[\cite{FF:Invariance}, Theorem 4.7]  
\label{local-to-global}
A holomorphic map $\pi:E\to B$ between complex manifolds is Oka if every point of $B$ has an open neighbourhood $U$ such that the restriction $\pi:\pi^{-1}(U)\to U$ is Oka.
\end{theorem}

It follows immediately that a holomorphic fibre bundle whose fibres are Oka manifolds is an Oka map.  The next result gives the weakest currently known geometric sufficient conditions for a holomorphic map to be Oka.  It is the strongest available statement of Gromov's Oka principle.

\begin{theorem}[\cite{FF:OkaMaps}, Corollary 1.2]  
\label{strongest-theorem}
{\rm(i)}  A stratified holomorphic fibre bundle with Oka fibres has the parametric Oka property.

{\rm (ii)}  A stratified subelliptic submersion has the parametric Oka property.

Thus, a holomorphic submersion of one of these two types is Oka if and only if it is a topological fibration.
\end{theorem}

\begin{example}
\label{eremenko}
There exists a holomorphic submersion $\pi:E\to B$ that is a smooth fibre bundle such that the fibre $\pi^{-1}(b)$ is an Oka manifold for every $b\in B$,
but $\pi$ is not an Oka map.  Here is a simple example.  Let $g:\mathbb D\to \C$ be a smooth function.  Let $\pi:E_g=\mathbb D\times \C\setminus\Gamma_g \to \mathbb D$ be the projection, where $\Gamma_g$ denotes the graph of $g$.  Clearly, $\pi$ is smoothly trivial and each fibre $\pi^{-1}(z)\cong\C\setminus\{g(z)\} \cong \C^*$ is an Oka manifold.  However, {\it if $\pi$ is an Oka map, then $g$ is holomorphic}.  Indeed, if $\pi$ is Oka, then the smooth lifting $f:\mathbb D\times\C^*\to E_g$, $(z,w)\mapsto (z,w+g(z))$, of the projection $p:\mathbb D\times\C^*\to\mathbb D$ by $\pi$ can be deformed to a holomorphic lifting $h: \mathbb D\times \C^*\to E_g$.
\[ \xymatrix{ & E_g \ar[d]^{\pi} \\ 
	            \mathbb D\times \C^* \ar@<1ex>[ur]^f \ar[ur]_h \ar[r]_{\ \ \ p} & \mathbb D } \]
For each $z\in\mathbb D$, $g(z)$ is the missing value in the range of the holomorphic map $h(z,\cdot):\C^*\to \C$. A theorem of A.~Eremenko \cite{Eremenko} now implies that $g$ is holomorphic.
\end{example}


\section{The homotopy-theoretic viewpoint}
\label{homotopy}

\noindent
Gromov's Oka principle has a natural home in structures provided by abstract homotopy theory.  This point of view was developed in the papers \cite{Larusson1, Larusson2, Larusson3, Larusson5}.  The Oka property of holomorphic maps (and in particular of complex manifolds) turns out to have a natural and rigorous homotopy-theoretic interpretation. 

Abstract homotopy theory, also known as homotopical algebra, was founded by D.\ Quillen in his 1967 monograph \cite{Quillen1967}.  The fundamental notion of the theory is the concept of a model category, or a model structure on a category.  Roughly speaking, a model structure is an abstraction of the essential features of the category of topological spaces that make ordinary homotopy theory possible.  We can only very briefly review the basic notions of abstract homotopy theory that are relevant here.  For more information we refer the reader to \cite{Dwyer-Spalinski, Hovey}.

A model category is a category with all small limits and colimits and three distinguished classes of maps, called weak equivalences or acyclic maps, fibrations, and cofibrations, such that the following axioms hold.
\begin{itemize}
\item  If $f$ and $g$ are composable maps, and two of $f$, $g$, $f\circ g$ are acyclic, so is the third.
\item  The classes of weak equivalences, fibrations, and cofibrations are closed under retraction.  (Also, it follows from the axioms that the composition of fibrations is a fibration, and the pullback of a fibration by an arbitrary map is a fibration.) 
\item A lifting $B\to X$ exists in every square
\[\xymatrix{ A \ar[r] \ar[d]_{\text{cofibration}} & X \ar[d]^{\text{fibration}} \\ B \ar[r] \ar@{-->}[ur] & Y }\]
if one of the vertical maps is acyclic.
\item  Every map can be functorially factored as
\[\text{acyclic fibration}\,\circ\,\text{cofibration}\] 
and as
\[\text{fibration}\,\circ\,\text{acyclic cofibration.}\]
\end{itemize}

There are many examples of model categories.  A fundamental example, closely related to the category of topological spaces, is the category of simplicial sets.  Simplicial sets are combinatorial objects that have a homotopy theory equivalent to that of topological spaces, but tend to be more useful or at least more convenient than topological spaces for various homotopy-theoretic purposes.  In homotopy-theoretic parlance, the distinction between topological spaces and simplical sets is blurred and the latter are often referred to as spaces.  For an introduction to simplicial sets, we refer the reader to \cite{Goerss-Jardine, May}.

The following key observations vastly expand the scope of abstract homotopy theory.
\begin{itemize}
\item  Not only can we do homotopy theory with individual spaces, but also with diagrams or sheaves of them.
\item  Manifolds and varieties can be thought of as sheaves of spaces, so we can do homotopy theory with them too (the general idea is known as the Yoneda lemma).
\end{itemize}
This line of thought has found a spectacular application in V.\ Voevodsky's homotopy theory of schemes and the resulting proof of the Milnor conjecture \cite{Voevodsky}.

We wish to embed the category $\sM$ of complex manifolds and holomorphic maps into a suitable model category.  Every complex manifold $X$ defines a presheaf $\sO(\cdot, X)$ on the full subcategory $\sS$ of Stein manifolds.  The presheaf consists of the set $\sO(S, X)$ of holomorphic maps $S\to X$ for each Stein manifold $S$, along with the precomposition map $\sO(S_2, X)\to\sO(S_1, X)$ induced by each holomorphic map $S_1\to S_2$ between Stein manifolds.  The presheaf $\sO(\cdot, X)$ determines $X$, so we have an embedding, in fact a full embedding, of $\sM$ into the category of presheaves of sets on $\sS$.

Each set $\sO(S, X)$ carries the compact-open topology.  A map between such sets defined by pre- or postcomposition by a holomorphic map is continuous.  We may therefore consider a complex manifold $X$ as a presheaf of topological spaces on $\sS$.  This presheaf has the property that as a holomorphic map $S_1\to S_2$ between Stein manifolds is varied continuously in $\sO(S_1,S_2)$, the induced precomposition map $\sO(S_2, X)\to\sO(S_1, X)$ varies continuously as well.  We would like to do homotopy theory with complex manifolds viewed as presheaves with this property.

Somewhat unexpectedly, as explained in \cite{Larusson2}, \S 3, there are solid reasons, beyond mere convenience, to rephrase the above entirely in terms of simplicial sets.  For the technical terms that follow, we refer the reader to \cite{Larusson2} and the references cited there.  To summarize, we turn $\sS$ into a simplicial site and obtain an embedding of $\sM$ into the category $\mathfrak S$ of prestacks on $\sS$.  The basic homotopy theory of prestacks on a simplicial site was developed by B.\ To\"en and G.\ Vezzosi for use in algebraic geometry \cite{Toen-Vezzosi}.  A new model structure on $\mathfrak S$, called the {\it intermediate structure} and based on ideas of J.\ F.\ Jardine, later published in \cite{Jardine}, was constructed in \cite{Larusson2}.  It is in this model structure that Gromov's Oka principle finds a natural home.

The main results of \cite{Larusson2} along with Theorem 6 of \cite{Larusson3} can be summarized as follows.

\begin{theorem}
The category of complex manifolds and holomorphic maps can be embedded into a model category such that:
\begin{itemize}
\item  a holomorphic map is acyclic when viewed as a map in the ambient model category if and only if it is a homotopy equivalence in the usual topological sense.
\item  a holomorphic map is a fibration if and only if it is an Oka map.  In particular, a complex manifold is fibrant if and only if it is Oka.
\item  a complex manifold is cofibrant if and only if it is Stein.
\item  a Stein inclusion is a cofibration.
\end{itemize}
\end{theorem}

A characterization of those holomorphic maps that are cofibrations is missing from this result.  It may be that Stein inclusions and biholomorphisms are the only ones.

Knowing that Oka maps are fibrations in a model structure helps us understand and predict their behaviour.  For example, by abstract nonsense, in any model category, the source of a fibration with a fibrant target is fibrant.  It follows that the source of an Oka map with an Oka target is Oka.  On the other hand, the fact that the image of an Oka map with an Oka source is Oka is a somewhat surprising feature of Oka theory not predicted by abstract nonsense, the reason being that the Oka property can be detected using Stein inclusions of the special kind $T\hookrightarrow \C^n$ with $T$ contractible.

Further connections between Oka theory and homotopy theory were pursued in \cite{Larusson4}.  The singular set $sX$ of a topological space $X$ is a simplicial set whose $n$-simplices for each $n\geq 0$ are the continuous maps into $X$ from the standard $n$-simplex
\[T_n=\{(t_0,\dots,t_n)\in\R^{n+1}:t_0+\dots+t_n=1,\, t_0,\dots,t_n\geq 0\}.\]
The singular set carries the homotopy type of $X$.

The affine singular set $eX$ of a complex manifold $X$ was defined in \cite{Larusson4} as a simplicial set whose $n$-simplices for each $n\geq 0$ are the holomorphic maps into $X$ from the affine $n$-simplex 
\[A_n=\{(t_0,\dots,t_n)\in\C^{n+1}:t_0+\dots+t_n=1\},\]
viewed as a complex manifold biholomorphic to $\C^n$.  If $X$ is Brody hyperbolic, then $eX$ is discrete and carries no topological information about $X$.

A holomorphic map $A_n\to X$ is determined by its restriction to $T_n\subset A_n$, so we have a monomorphism, that is, a cofibration $eX\hookrightarrow sX$ of simplical sets.  When $X$ is Oka, $eX$, which is of course much smaller than $sX$, carries the homotopy type of $X$.  More precisely, the cofibration $eX\hookrightarrow sX$ is the inclusion of a strong deformation retract (\cite{Larusson4}, Theorem 1).  Even for complex Lie groups, this result appears not to have been previously known.


\section{Open problems}
\label{open}

\noindent
Oka theory is a young field with many basic problems still open.  A major theme is to clarify the boundaries of the class of Oka manifolds.  We start with two frustratingly simple-looking questions.

\noindent
{\bf A.\ }  Let $B$ be a closed ball in $\C^n$, $n\geq 2$ (or more generally a compact convex set).  Is the complement $\C^n\setminus B$ Oka?  What makes this problem particularly intriguing is the absence of any obvious obstructions.  Indeed, $\C^n\setminus B$ is a union of Fatou-Bieberbach domains \cite{Rosay-Rudin}.

\noindent
{\bf B.\ }  Is the complement of a smooth cubic curve in $\P^2$ Oka?  The complement is known to be dominable by $\C^2$ \cite{Buzzard-Lu}.

\smallskip
The following are broader questions in a similar vein.

\noindent
{\bf C.\ }  Which complex surfaces of non-general type are Oka?  In particular, how about Kummer surfaces and, more generally, K3 surfaces?  Kummer surfaces are dominable by $\C^2$ and they are dense in the moduli space of all K3 surfaces.

A Kummer surface $X$ admits a finite branched covering $Y\to X$ where $Y$ is a complex two-dimensional torus blown up at 16 points.  The universal covering space $\widetilde Y$ of $Y$ is $\C^2$ (the universal covering of the torus) blown up along a tame discrete sequence (the preimage points of the 16 points in $X$), so $\widetilde Y$ is weakly subelliptic.  Hence, $\widetilde Y$ is Oka, so  by Corollary \ref{covers-quotients}, $Y$ is also Oka.  It remains to be seen whether the Oka property passes down from $Y$ to its ramified quotient space $X$.

\noindent
{\bf D.\ }  If $X$ is a complex manifold, $p\in X$, and $X\setminus\{p\}$ is Oka, is $X$ Oka?  Conversely, if $X$ is Oka and $\dim X\geq 2$, is $X\setminus\{p\}$ Oka?

\smallskip
Continuing in the same vein, we inquire about possible counterparts in Oka theory to well-known statements in hyperbolic geometry.

\noindent
{\bf E.\ }  S.\ Kobayashi conjectured that if $X$ is a very general hypersurface in $\P^n$ of sufficiently high degree, then $X$ and $\P^n\setminus X$ are hyperbolic.  There has been considerable progress on this conjecture, which we shall not review here, but as far as we know a complete proof has not appeared.  Is a smooth hypersurface in $\P^n$ of sufficiently low degree (greater than~$1$) Oka?  Is the complement of a generic hypersurface in $\P^n$ of sufficiently low degree Oka?

\noindent
{\bf F.\ }  It is known that hyperbolicity of compact complex manifolds is stable under small deformations (\cite{Kobayashi}, Theorem 3.11.1).  How does the Oka property of compact complex manifolds behave with respect to deformations?

\smallskip
The next two problems concern a geometric characterization of the Oka property (for manifolds for now; later, one would want to generalize to maps).  The essence of Gromov's Oka principle is the implication elliptic $\Rightarrow$ Oka.  Gromov proved the converse for Stein manifolds (\cite{Gromov:OP}, 3.2.A; see also \cite{Larusson3}, Theorem 2).  For manifolds in general, this is a fundamental open question.

\noindent
{\bf G.\ }  Is every Oka manifold elliptic (or subelliptic or weakly subelliptic)?

Gromov's result that ellipticity is equivalent to the Oka property for Stein manifolds has been generalized to a much larger class of manifolds using a geometric structure somewhat more involved than a dominating spray (\cite{Larusson3}; see also \cite{Larusson5}, \S 3).  

A complex manifold is said to be {\it good} if it is the image of an Oka map from a Stein manifold, and {\it very good} if it carries a holomorphic affine bundle whose total space is Stein (these definitions are slightly different from those in \cite{Larusson3}).

The simplest examples of good manifolds that are not Stein are the projective spaces.  Namely, let $Q_n$ be the complement in $\P^n\times\P^n$ of the hypersurface
\[\big\{([z_0,\dots,z_n], [w_0,\dots,w_n]):z_0 w_0+\dots+z_n w_n=0\big\}.\]
This hypersurface is the preimage of a hyperplane by the Segre embedding $\P^n\times \P^n\to \P^{n^2+2n}$, so $Q_n$ is Stein.  Let $\pi$ be the projection $Q_n\to\P^n$ onto the first component.  It is easily seen that $\pi$ has the structure of a holomorphic affine bundle with fibre $\C^n$.  Thus, $\P^n$ is very good.  (This observation is called the Jouanolou trick in algebraic geometry.)

The class of good manifolds contains all Stein manifolds and all quasi-projective manifolds and is closed under taking submanifolds, products, covering spaces, finite branched covering spaces, and complements of analytic hypersurfaces.  The same is true of the class of very good manifolds.

A good manifold is Oka if and only if it is the image of an Oka map from an elliptic manifold.  A very good manifold is Oka if and only if it carries an affine bundle whose total space is elliptic: this is a purely geometric characterization of the Oka property that holds, for example, for all quasi-projective manifolds.

\noindent
{\bf H.\ }  Is every complex manifold good, or even very good?

\smallskip
We conclude with an assortment of further problems.

\noindent
{\bf I.\ }   Are the affine spaces $\C^n$ the only contractible Stein Oka manifolds?  S.\ Kaliman and F.\ Kutzschebauch have produced many smooth affine algebraic varieties that are diffeomorphic to affine space and have the algebraic density property, so they are Oka \cite{Kaliman-Kutzschebauch}.  It is not known whether these manifolds are biholomorphic to affine space.

\noindent
{\bf J.\ }   If we restrict the convex compact subset $K$ in CAP to be a ball, is the resulting property equivalent to CAP?  This question was posed by Gromov (\cite{Gromov:OP}, p.\ 881). 

\noindent
{\bf K.\ }  Let $E\to B$ be a holomorphic fibre bundle with fibre $Y$.  If $E$ is Oka, does it follow that $B$ and $Y$ are Oka as well?  It would be surprising if the answer was affirmative, but no counterexample is known.

\noindent
{\bf L.\ }  Is there a reasonable way to extend the Oka property from complex manifolds to reduced complex spaces?  The authors have tried and found it problematic.

\noindent
{\bf M.\ }   Are there any restrictions on the homotopy type of a compact Oka manifold?  It would be amusing if simply connected compact Oka manifolds turned out to have elliptic homotopy type in the sense of rational homotopy theory.  This would imply, for example, that K3 surfaces are not Oka \cite{Amoros-Biswas}.


\begin{thebibliography}{99}

\bibitem{Amoros-Biswas}
{\scshape Amor\'os, J.; Biswas, I.}
Compact K\"ahler manifolds with elliptic homotopy type.  {\em Adv.\ Math.} \textbf{224} (2010), 1167--1182.  MR2628808, Zbl 1198.32008.

\bibitem{Arnold}
{\scshape Arnold, J.\ E., Jr.}
Local to global theorems in the theory of Hurewicz fibrations.  {\em Trans.\ Amer.\ Math.\ Soc.} \textbf{164} (1972), 179--188.  MR0295349 (45 \#4415), Zbl 0227.55014.

\bibitem{Brody}
{\scshape Brody, R.}
Compact manifolds and hyperbolicity.
{\em Trans.\ Amer.\ Math.\ Soc.} \textbf{235} (1978), 213--219.  MR0470252 (57 \#10010), Zbl 0416.32013.

\bibitem{Buzzard-Lu}
{\scshape Buzzard, G.; Lu, S.\ S.\ Y.}
Algebraic surfaces holomorphically dominable by $\C^2$. {\em Invent. Math.} \textbf{139} (2000), 617--659.  MR1738063 (2000m:14040), Zbl 0967.14025.

\bibitem{Dwyer-Spalinski}
{\scshape Dwyer, W.\ G.; Spali\'nski, J.}
Homotopy theories and model categories. 
{\em Handbook of algebraic topology.}  {\em North-Holland,} 1995, 73--126.  MR1361887 (96h:55014), Zbl 0869.55018.

\bibitem{Eliashberg1}
{\scshape Eliashberg, Y.}
Topological characterization of Stein manifolds of dimension $>2$.
{\em Internat.\ J.\ Math.} \textbf{1} (1990), 29--46.  MR1044658 (91k:32012), Zbl 0699.58002.

\bibitem{Eliashberg-Gromov1}
{\scshape Eliashberg, Y.; Gromov, M.}
Nonsingular mappings of Stein manifolds.
{\em Funkcional.\ Anal.\ i Prilo\v zen.} \textbf{5} (1971), 82--83.  MR0301236 (46 \#394), Zbl 0234.32011.

\bibitem{Eliashberg-Gromov2}
{\scshape Eliashberg, Y.; Gromov, M.}
Embeddings of Stein manifolds of dimension $n$ into the affine space of dimension $3n/2+1$. 
{\em Ann.\ Math.} (2) \textbf{136} (1992), 123--135.  MR1173927 (93g:32037), Zbl 0758.32012.

\bibitem{Eremenko}
{\scshape Eremenko, A.}
Exceptional values in holomorphic families of entire functions.
{\em Michigan Math.\ J.} \textbf{54} (2006), 687--696.  MR2280501 (2008b:32003), Zbl 1155.32005.

\bibitem{FF:Actions}
{\scshape Forstneri\v c, F.}
Actions of $(\mathbf R,+)$ and $(\mathbf C,+)$ on complex manifolds.
{\em Math.\ Z.} \textbf{223} (1996), 123--153.  MR1408866 (97i:32041), Zbl 0872.32021.

\bibitem{FF:Intersections}
{\scshape Forstneri\v c, F.}
On complete intersections.  {\em Ann.\ Inst.\ Fourier (Grenoble)} \textbf{51} (2001), 497--512.  MR1824962 (2002c:32044), Zbl 0991.32008.

\bibitem{FF:Subelliptic}
{\scshape Forstneri\v c, F.}
The Oka principle for sections of subelliptic submersions.  {\em Math.\ Z.} \textbf{241} (2002), 527--551.  MR1938703 (2003i:32043), Zbl 1023.32008.

\bibitem{FF:Noncritical}
{\scshape Forstneri\v c, F.}
Noncritical holomorphic functions on Stein manifolds.  {\em Acta Math.} \textbf{191} (2003), 143--189.  MR2051397 (2005b:32021), Zbl 1064.32021.

\bibitem{FF:EOP}
{\scshape Forstneri\v c, F.}
Extending holomorphic mappings from subvarieties in Stein manifolds.
{\em Ann.\ Inst.\ Fourier (Grenoble)} \textbf{55} (2005), 733--751.  MR2149401 (2006c:32012), Zbl 1076.32003.

\bibitem{FF:Flexibility}
{\scshape Forstneri\v c, F.}
Holomorphic flexibility properties of complex manifolds.
{\em Amer.\ J.\ Math.} \textbf{128} (2006), 239--270.  MR2197073 (2006k:32024), Zbl 1171.32303.

\bibitem{FF:CAP}
{\scshape Forstneri\v c, F.}
Runge approximation on convex sets implies the Oka property.
{\em Ann.\ Math.} (2) \textbf{163} (2006), 689--707.  MR2199229 (2006j:32011), Zbl 1103.32004.

\bibitem{FF:OkaManifolds} 
{\scshape Forstneri\v c, F.}
Oka manifolds.  {\em C.\ R.\ Acad.\ Sci.\ Paris, Ser.\ I} \textbf{347} (2009), 1017--1020.  MR2554568, Zbl 1175.32005.

\bibitem{FF:OkaMaps}
{\scshape Forstneri\v c, F.}
Oka maps.  {\em C.\ R.\ Acad.\ Sci.\ Paris, Ser.\ I} \textbf{348} (2010), 145--148.  MR2600066, Zbl pre05682595.

\bibitem{FF:Kohn}
{\scshape Forstneri\v c, F.}
The Oka principle for sections of stratified fiber bundles.
{\em Pure Appl.\ Math.\ Quarterly} \textbf{6} (2010), 843--874.  MR2677316.

\bibitem{FF:Invariance}
{\scshape Forstneri\v c, F.}
Invariance of the parametric Oka property.
{\em Complex analysis.}  Trends in Mathematics.  {\em Birkh\"auser,} 2010, 125--144.  Zbl pre05762421.

\bibitem{Forstneric-Prezelj2000}
{\scshape Forstneri\v c, F.; Prezelj, J.}
Oka's principle for holomorphic fiber bundles with sprays.  {\em Math.\ Ann.} \textbf{317} (2000), 117--154.  MR1760671 (2001m:32040), Zbl 0964.32017.

\bibitem{Forstneric-Prezelj2001}
{\scshape Forstneri\v c, F.; Prezelj, J.}
Extending holomorphic sections from complex subvarieties.  {\em Math.\ Z.} \textbf{236} (2001), 43--68.  MR1812449 (2002b:32017), Zbl 0968.32005.

\bibitem{Forstneric-Prezelj2002}
{\scshape Forstneri\v c, F.; Prezelj, J.} 
Oka's principle for holomorphic submersions with sprays.  {\em Math.\ Ann.} \textbf{322} (2002), 633--666.  MR1905108 (2003d:32027), Zbl 1011.32006.

\bibitem{Forstneric-Slapar2007a}
{\scshape Forstneri\v c, F.; Slapar, M.}
Stein structures and holomorphic mappings.  {\em Math.\ Z.} \textbf{256} (2007), 615--646.  MR2299574 (2008g:32036), Zbl 1129.32013.

\bibitem{Forstneric-Slapar2007b}
{\scshape Forstneri\v c, F.; Slapar, M.}
Deformations of Stein structures and extensions of holomorphic mappings.  {\em Math.\ Res.\ Lett.} \textbf{14} (2007), 343--357.  MR2318630 (2008d:32023), Zbl 1134.32007.

\bibitem{Forstneric-Wold1}
{\scshape Forstneri\v c, F.; Wold, E.\ F.}
Bordered Riemann surfaces in $\C^2$. 
{\em J.\ Math.\ Pures Appl.} \textbf{91} (2009), 100--114.  MR2487902 (2010b:32008), Zbl 1157.32010.

\bibitem{Forstneric-Wold2}
{\scshape Forstneri\v c, F.; Wold, E.\ F.}
Fibrations and Stein neighborhoods.
{\em Proc.\ Amer.\ Math.\ Soc.} \textbf{138} (2010), 2037--2042.  MR2596039 (2011a:32016), Zbl 1192.32008.

\bibitem{Goerss-Jardine}
{\scshape Goerss P.\ G.; Jardine, J.\ F.}
Simplicial homotopy theory.  Progress in Mathematics, 174.  {\em Birkh\"auser,} 1999.  MR1711612 (2001d:55012), Zbl 0949.55001.

\bibitem{Gompf1}
{\scshape Gompf, R.\ E.}
Handlebody construction of Stein surfaces.
{\em Ann.\ Math.} (2) \textbf{148} (1998), 619--693.  MR1668563 (2000a:57070), Zbl 0919.57012.

\bibitem{Gompf2}
{\scshape Gompf, R.\ E.}
Stein surfaces as open subsets of $\C^2$.
{\em J.\ Symplectic Geom.} \textbf{3} (2005), 565--587.  MR2235855 (2007f:32031), Zbl 1118.32011.

\bibitem{Grauert3}
{\scshape Grauert, H.}
Analytische Faserungen \"uber holomorph-vollst\"an\-digen R\"au\-men.
{\em Math.\ Ann.} \textbf{135} (1958), 263--273.  MR0098199 (20 \#4661), Zbl 0081.07401. 

\bibitem{Grauert-Remmert}
{\scshape Grauert, H.; Remmert, R.}
Theorie der Steinschen R\"aume.  Grundlehren der mathematischen Wissenschaften, 227.  {\em Springer-Verlag,} 1977.  Translation: Theory of Stein spaces.  Grundlehren der mathematischen Wissenschaften, 236. {\em Springer-Verlag,} 1979.  MR0513229 (80j:3200), Zbl 0379.32001.

\bibitem{Gromov:PDR}
{\scshape Gromov, M.}
Partial differential relations.  Ergebnisse der Mathematik und ihrer Grenzgebiete (3), 9. {\em Springer-Verlag,} 1986.  MR0864505 (90a:58201), Zbl 0651.53001.

\bibitem{Gromov:OP}
{\scshape Gromov, M.}
Oka's principle for holomorphic sections of elliptic bundles.
{\em J.\ Amer.\ Math.\ Soc.} \textbf{2} (1989), 851-897.  MR1001851 (90g:32017), Zbl 0686.32012.

\bibitem{Gunning-Rossi}
{\scshape Gunning, R.\ C.; Rossi, H.}
Analytic functions of several complex variables.  {\em Prentice-Hall,} 1965.  MR0180696 (31 \#4927), Zbl 0141.08601.

\bibitem{Henkin-Leiterer:Oka}
{\scshape Henkin, G.\ M.; Leiterer, J.}
The Oka-Grauert principle without induction over the base dimension.
{\em Math.\ Ann.} \textbf{311} (1998), 71--93.  MR1624267 (99f:32048), Zbl 0955.32019.

\bibitem{Hormander-SCV}
{\scshape H\"ormander, L.}
An introduction to complex analysis in several variables.  Third edition.
North-Holland Mathematical Library, 7. {\em North Holland,} 1990.  MR1045639 (91a:32001), Zbl 0685.32001.

\bibitem{Hovey}
{\scshape Hovey, M.}
Model categories. Mathematical Surveys and Monographs, 63. {\em Amer.\ Math.\ Soc.,} 1999.  MR1650134 (99h:55031), Zbl 0909.55001.

\bibitem{Ivarsson-Kutzschebauch}
{\scshape Ivarsson, B.; Kutzschebauch, F.}
A solution of Gromov's Vaserstein problem.
{\em C.\ R.\ Acad.\ Sci.\ Paris, Ser.\ I} \textbf{346} (2008), 1239--1243. MR2473300 (2010a:32022), Zbl 1160.32017.

\bibitem{Jardine}
{\scshape Jardine, J.\ F.}
Intermediate model structures for simplicial presheaves.  {\em Canad.\ Math.\ Bull.} \textbf{49} (2006), 407--413.  MR2252262 (2007d:18021), Zbl 1107.18007.

\bibitem{Kaliman-Kutzschebauch}
{\scshape Kaliman, S.; Kutzschebauch, F.}
Density property for hypersurfaces $UV=P(\overline X)$.  {\em Math.\ Z.} \textbf{258} (2008), 115--131.  MR2350038 (2008k:32062), Zbl 1133.32012.

\bibitem{Kaliman-Zaidenberg}
{\scshape Kaliman, S.; Zaidenberg, M.}
A tranversality theorem for holomorphic mappings and stability of Eisenman-Kobayashi measures.  
{\em Trans.\ Amer.\ Math.\ Soc.} \textbf{348} (1996), 661--672.  MR1321580 (96g:32043), Zbl 0851.32018.

\bibitem{Kobayashi}
{\scshape Kobayashi, S.}
Hyperbolic complex spaces. Grundlehren der mathematischen Wissenschaften, 318. {\em Springer-Verlag,} 1998.  MR1635983 (99m:32026), Zbl 0917.32019.

\bibitem{Kobayashi-Ochiai}
{\scshape Kobayashi, S.; Ochiai, T.}
Meromorphic mappings onto compact complex spaces of general type.
{\em Invent. Math.} \textbf{31} (1975), 7--16.  MR0402127 (53 \#5948), Zbl 0331.32020.

\bibitem{Kolaric}
{\scshape Kolari\v c, D.}
Parametric h-principle for holomorphic immersions with approximation.  Preprint, {\tt arXiv:1005.1274} (2010).

\bibitem{Larusson1}
{\scshape L\'arusson, F.}
Excision for simplicial sheaves on the Stein site and Gromov's Oka principle.
{\em Internat.\ J.\ Math.} \textbf{14} (2003), 191--209.  MR1966772 (2004b:32041), Zbl 1078.32017.

\bibitem{Larusson2}
{\scshape L\'arusson, F.}
Model structures and the Oka principle.
{\em J.\ Pure Appl.\ Algebra} \textbf{192} (2004), 203--223.  MR2067196 (2005e:32043), Zbl 1052.32020.

\bibitem{Larusson3}
{\scshape L\'arusson, F.}
Mapping cylinders and the Oka principle.
{\em Indiana Univ.\ Math.\ J.} \textbf{54} (2005), 1145--1159.  MR2164421 (2006f:32035), Zbl 1085.32011.

\bibitem{Larusson4}
{\scshape L\'arusson, F.}
Affine simplices in Oka manifolds.
{\em Documenta Math.} \textbf{14} (2009), 691--697.  MR2578806 (2010m:32029), Zbl 1200.32016.

\bibitem{Larusson5}
{\scshape L\'arusson, F.}
Applications of a parametric Oka principle for liftings.
{\em Complex analysis.}  Trends in Mathematics.  {\em Birkh\"auser,} 2010, 205--211.  Zbl pre05762426.

\bibitem{Larusson6}
{\scshape L\'arusson, F.}
What is an Oka manifold?
{\em Notices Amer.\ Math.\ Soc.} \textbf{57} (2010), 50--52.  MR2590114, Zbl 1191.32004.

\bibitem{Majcen-2009}
{\scshape Majcen, I.}
Embedding certain infinitely connected subsets of bordered Riemann surfaces properly into $\C^2$.  {\em J.\ Geom.\ Anal.} \textbf{19} (2009), 695--707.  MR2496573 (2010j:32023), Zbl 1172.32300.

\bibitem{May}
{\scshape May, J.\ P.}
Simplicial objects in algebraic topology.  Van Nostrand Mathematical Studies, 11.  {\em D.\ Van Nostrand Co.,} 1967.  Chicago Lectures in Mathematics.  {\em University of Chicago Press,} 1992.  MR0222892 (36 \#5942), Zbl 0165.26004.

\bibitem{Oka1939}
{\scshape Oka, K.}
Sur les fonctions des plusieurs variables.  III: Deuxi\`eme probl\`eme de Cousin.
{\em J.\ Sc.\ Hiroshima Univ.} \textbf{9} (1939), 7--19.

\bibitem{Prezelj2010}
{\scshape Prezelj, J.}
A relative Oka-Grauert principle for holomorphic submersions over 1-convex spaces.  {\em Trans.\ Amer.\ Math.\ Soc.} \textbf{362} (2010), 4213--4228.  MR2608403, Zbl pre05770877.

\bibitem{Quillen1967}
{\scshape Quillen, D.}
Homotopical algebra.  Lecture Notes in Mathematics, 43.  {\em Springer-Verlag,} 1967.  MR0223432 (36 \#6480), Zbl 0168.20903.

\bibitem{Rosay-Rudin}
{\scshape Rosay, J.-P.; Rudin, W.}
Holomorphic maps from $\C^n$ to $\C^n$. {\em Trans. Amer. Math. Soc.} \textbf{310} (1988), 47--86.  MR0929658 (89d:32058), Zbl 0708.58003.

\bibitem{Schurmann}
{\scshape Sch\"urmann, J.}
Embeddings of Stein spaces into affine spaces of minimal dimension.
{\em Math.\ Ann.} \textbf{307} (1997), 381--399.  MR1437045 (98a:32011), Zbl 0881.32007.

\bibitem{Siu1976}
{\scshape Siu, Y.-T.}
Every Stein subvariety admits a Stein neighborhood.
{\em Invent.\ Math.} \textbf{38} (1976), 89--100.  MR0435447 (55 \#8407), Zbl 0343.32014.

\bibitem{Toen-Vezzosi}
{\scshape To\"en, B.; Vezzosi, G.}
Homotopical algebraic geometry. I. Topos theory. {\em Adv.\ Math.} \textbf{193} (2005), 257--372.  MR2137288 (2007b:14038), Zbl 1120.14012.

\bibitem{Voevodsky}
{\scshape Voevodsky, V.}
$\mathbf A^1$-homotopy theory.  Proceedings of the International Congress of Mathematicians, vol.\ I (Berlin, 1998). {\em Documenta Math.}, extra vol.\ I (1998), 579--604.  MR1648048 (99j:14018), Zbl 0907.19002.

\end{thebibliography}
\end{document}